\documentclass[final]{siamltex}

\usepackage{amsmath,amssymb,epsfig}
\title{Quasi-optimal robust stabilization of control systems}

\author{Christophe Prieur\thanks{LAAS-CNRS, 7 avenue du Colonel Roche
31077 Toulouse, France ({\tt christophe.prieur@laas.fr}).}
\and Emmanuel Tr\'elat\thanks{Univ.\ Paris-Sud, Labo.\ Math., UMR
8628, Bat.\ 425, 91405 Orsay Cedex, France
({\tt emmanuel.trelat@math.u-psud.fr}).}}

\def\5{7} 
\def\4{6} 
\def\R{\R}

\def\NNN{\mathcal{N}}
\def\R{\rm{I\kern-0.21emR}}
\def\N{\textrm{I\kern-0.21emN}}

\def\Om{\Omega}
\def\al{\alpha}

\def\FF{\mathcal{F}}
\def\clo{\mathrm{clos}}

\def\infini{\omega}

\def\ms{\mapsto}
\def\fa{\forall}

\newcommand{\eps}{\varepsilon}

\def\sup{\mathrm{sup}}

\def\overm{\overline{m}}
\def\ra{\rightarrow}

\newtheorem{remark}[theorem]{Remark}

\renewcommand{\geq}{\geqslant}
\renewcommand{\leq}{\leqslant}
\renewcommand{\rm}{\textrm}
\renewcommand{\it}{\textit}

\begin{document}
\maketitle

\begin{abstract}
In this paper, we investigate the problem of semi-global minimal
time robust stabilization of analytic control systems with
controls entering linearly, by means of a hybrid state feedback
law. It is shown that, in the absence of minimal time singular
trajectories, the solutions of the closed-loop system converge to
the origin in quasi minimal time (for a given bound on the
controller) with a robustness property with respect to small
measurement noise, external disturbances and actuator noise.
\end{abstract}



\pagestyle{myheadings}
\thispagestyle{plain}
\markboth{C.~PRIEUR AND E.~TRELAT}{Quasi-optimal robust stabilization}

\section{Introduction}\label{intro}
Let $m$ and $n$ be two positive integers.
Consider on $\R^n$ the control system
\begin{equation}\label{systeme}
\dot{x}(t) = \sum_{i=1}^m u_i(t) f_i(x(t)),
\end{equation}
where $f_1,\ldots,f_m$ are analytic vector fields on $\R^n$, and where the
control function
$u(\cdot)=(u_1(\cdot),\ldots,u_m(\cdot))$ satisfies the constraint
\begin{equation}\label{contrainte}
\sum_{i=1}^m u_i(t)^2\leq 1.
\end{equation}
All results of this paper still hold on a Riemannian analytic manifold of
dimension $n$, which is connected and complete. However, for the sake of
simplicity, our results are stated in $\R^n$. Let $\bar{x}\in \R^n$.
The system (\ref{systeme}),
together with the constraint (\ref{contrainte}), is said
\textit{globally asymptotically stabilizable} at the point
$\bar{x}$, if, for each point $x\in\R^n$,
there exists a control law satisfying the constraint (\ref{contrainte}) such
that the solution of (\ref{systeme}) associated to this control law
and starting
from $x$ tends to $\bar{x}$ as $t$ tends to $+\infty$.

This asymptotic stabilization problem has a long history and has
been widely investigated. Note that, due to \textit{Brockett's
condition} \cite[Theorem 1, (iii)]{bro}, if $m<n$, and if the
vector fields $f_1,\ldots, f_m$ are independent, then there does
not exist any continuous stabilizing feedback law for
(\ref{systeme}). However several control laws have been derived
for such control systems (see for instance
\cite{ast_holo_EJC:98,kolm_mccl:95} and references therein).

The \textit{robust asymptotic stabilization problem}
is under current and
active research. Many notions of controllers have been introduced
to treat this problem, such as discontinuous sampling feedbacks
\cite{clrs,son:99}, time varying control laws
\cite{comu,coron92,mosa,mosa2}, patchy feedbacks (as in
\cite{anbr02}), SRS feedbacks \cite{riffordSRS}, ..., enjoying
different robustness properties depending on the errors under
consideration.

In the present paper, we consider feedback laws having both discrete and
continuous components, which generate
closed-loop systems with \textit{hybrid} terms (see for instance
\cite{bm,tav1}).
Such feedbacks appeared first in \cite{MCSS} to stabilize nonlinear
systems having a priori no discrete state. They consist in defining a
\textit{switching strategy} between several smooth control laws
defined on a partition of the state space. Many results on the
stabilization problem of nonlinear systems by means of hybrid controllers
have been recently established (see for instance
\cite{bra,piccoli2,go-te,nolcos04,liberzon,lygeros,piccoli3,ye}).
The notion of solution,
connected with the robustness problem, is by now well defined in the hybrid
context (see \cite{go-te,chained} among others).
Specific conditions for
the optimization can be found in the literature (see
\textit{e.g.}\ \cite{bac,piccoli1}).

The strategy of our paper is to combine a minimal time controller
that is smooth on a part of the state space, and other controllers
defined on the complement of this part, so as to provide a
\textit{quasi minimal time hybrid controller} by defining a
switching strategy between all control laws. The resulting hybrid
law enjoys a quasi minimal time property, and robustness with
respect to (small) measurement noise, actuator errors and external
disturbances.

More precisely, in a first step, we consider the minimal time
problem for the system (\ref{systeme}) with the constraint
(\ref{contrainte}), of steering a point $x\in \R^n$ to the point $\bar{x}$.
Note that this problem is solvable as soon as the \textit{Lie Algebra
  Rank Condition} holds for the $m$-tuple of vector fields
$(f_1,\ldots,f_m)$.
Of course, in general, it is impossible to compute explicitly the
minimal time feedback controllers for this problem. Moreover, Brockett's
condition implies that such control laws are not smooth whenever
$m<n$ and the vector fields $f_1,\ldots, f_m$ are independent.
This raises the problem of the regularity of optimal feedback laws.
The literature on this subject is immense.
In an analytic setting, the problem of determining the analytic
regularity of the minimal time function has been, among others,
investigated in \cite{Sussmann}. For systems of the
form (\ref{systeme}), it follows from \cite{Ag,AG,these} that the
minimal time function to $\bar{x}$ is \textit{subanalytic}, provided
there does not exist any
nontrivial singular minimal time trajectory starting from $\bar{x}$
(see \cite{Hardt,Hironaka} for a general definition of subanalytic
sets). This assumption holds generically for systems (\ref{systeme}),
whenever $m\geq 3$ (see \cite{CJT}).
In particular, this function is analytic outside
a stratified submanifold $\mathcal{S}$ of $\R^n$, of codimension
greater than or equal to $1$ (see \cite{Tamm}).
As a consequence, outside this submanifold it is possible to provide
an analytic minimal time feedback controller for the system
(\ref{systeme}), (\ref{contrainte}).
This optimal controller gives rise to trajectories never crossing
again the singular set $\mathcal{S}$.

Note that the analytic context is used so as to ensure stratification
properties, which do not hold a priori if the system is smooth only.
These properties are related to the notion of \textit{o-minimal
  category} (see \cite{DriesMiller}).

In a neighborhood of $\mathcal{S}$, we prove the existence of a set of
controllers steering the system (\ref{systeme}), (\ref{contrainte})
outside of this neighborhood in small time.

Then, in order to achieve a minimal time robust stabilization
procedure, using a hybrid feedback law, we define a suitable switching
strategy (more precisely, a hysteresis) between the minimal time
feedback controller and other controllers defined in a neighborhood of
$\mathcal{S}$. The resulting hybird system
has the following property: if the
state is close to the singular submanifold $\mathcal{S}$, the feedback
controller will push the state far enough from $\mathcal{S}$, in small
time; if the state is not too close to $\mathcal{S}$, then the
feedback controller will steer the system to $\bar{x}$ in minimal
time. Hence, the stabilization is quasi-optimal, and is proved
to enjoy robustness properties.

Note that we thus give an alternative solution, in the context of
hybrid systems using hysteresis, to a conjecture of
\cite[Conj.\ 1, p.\ 101]{Bressantorino}
concerning patchy feedbacks for smooth control
systems.\footnote{This conjecture on patchy feedbacks has been recently
considered in \cite{AnconaBressanPreprint}. In this preprint, written 
during the review process of the present work, the authors
prove, using a penalization method, a general result on stabilization by
means of patchy feedbacks of nonlinear
control systems in quasi-minimal time.}

\medskip

In a previous paper \cite{PrieurTrelat},
this program was achieved on the so-called Brockett system, for which
$n=3$, $m=2$, and, denoting $x=(x_1,x_2,x_3)$,
\begin{equation*}
f_1 = \frac{\partial}{\partial x_1}+x_2\frac{\partial}{\partial x_3},\ \
f_2 = \frac{\partial}{\partial x_2}-x_1\frac{\partial}{\partial x_3}.
\end{equation*}
In this case, there does not exist any nontrivial singular trajectory,
and the manifold $\mathcal{S}$ coincides with the axis $(0x_3)$. A simple
explicit hybrid strategy was described.
In contrast, in the present paper, we derive a general result that
requires a countable number of components in the definition of the
hysteresis hybrid feedback law.

\medskip

The paper is organized as follows.
In Section \ref{sec2}, we first recall some facts about the
minimal time problem for the system (\ref{systeme}) with
(\ref{contrainte}), and recall the definition of a singular
trajectory. Then, we give a notion of
solution adapted to hybrid feedback laws, and define the concept of
stabilization via a minimal time hybrid feedback law. The
main result, Theorem \ref{CH:TH1} in Section \ref{ch:main_sect},
states that, if there does not exist any
nontrivial singular minimal time trajectory of (\ref{systeme}),
(\ref{contrainte}), starting from $\bar{x}$,
then there exists a minimal time hybrid
feedback law stabilizing semi-globally the point $\bar{x}$ for the system
(\ref{systeme}), (\ref{contrainte}). Section \ref{secshortdescription}
describes the main ideas of the proof of the main result, and in
particular, contains two key lemmas. Section \ref{secproofmainthm} is
then devoted to the detailed proof of Theorem \ref{CH:TH1}, and
gathers all technical aspects needed to deal with hybrid systems: the
components of the hybrid feedback law, and a switching strategy
between both components are defined, and properties of the closed-loop
system are investigated.

The results in this work were announced in \cite{PTcdc05}.

\section{Definitions and main result}\label{sec2}
\subsection{The minimal time problem}
Consider the minimal time problem for the system (\ref{systeme}) with
the constraint (\ref{contrainte}).

Throughout the paper, we assume that the \textit{Lie Algebra Rank
  Condition} holds, that is, the Lie algebra spanned by the vector
fields $f_1,\ldots,f_m$ is equal to $\R^n$, at every point $x$ of
$\R^n$.

It is well known that, under this condition,
any two points of $\R^n$ can be joined by
a minimal time trajectory of (\ref{systeme}), (\ref{contrainte}).

Let $\bar{x}\in \R^n$ be fixed.
We denote by $T_{\bar{x}}(x)$ the minimal time needed to steer
the system (\ref{systeme}) with the constraint (\ref{contrainte}) from
a point $x\in\R^n$ to the point $\bar{x}$.

\begin{remark}
Obviously, the control function $u$ associated to a
minimal time
trajectory of (\ref{systeme}), (\ref{contrainte}), actually satisfies
$\sum_{i=1}^m u_i^2=1.$
\end{remark}

For $T>0$, let $\mathcal{U}_T$ denote the (open) subset of $u(\cdot)$ in
$L^\infty([0,T],\R^m)$
such that the solution of (\ref{systeme}), starting from $\bar{x}$
and associated to a control $u(\cdot) \in{\mathcal{U}}_T$, is
well defined on $[0,T]$. The mapping
$$\begin{array}{rcl}
 E_{\bar{x},T}:\ \ \ {\mathcal U}_T & \longrightarrow & \R^n \\
                    u(\cdot)  & \longmapsto & x(T),
\end{array}
$$
which to a control $u(\cdot)$ associates the end-point $x(T)$
of the corresponding solution $x(\cdot)$ of (\ref{systeme})
starting at $\bar{x}$, is called
{\textit{end-point mapping}} at the point $\bar{x}$, in time $T$; it is a
smooth mapping.

\begin{definition}\label{defsing}
A trajectory $x(\cdot)$ of (\ref{systeme}), so that $x(0)=\bar{x}$,
is said \textit{singular} on $[0,T]$ if
its associated control $u(\cdot)$ is a singular point of the end-point
mapping $E_{\bar{x},T}$ (\textit{i.e.}, if the Fr\'echet derivative of
$E_{\bar{x},T}$ at $u(\cdot)$ is not onto). The control $u(\cdot)$ is said
singular.
\end{definition}

\begin{remark}
If $x(\cdot)$ is singular on $[0,T]$, then it is singular on $[t_0,t_1]$,
for all $t_0,t_1\in [0,T]$ such that $t_0<t_1$.
\end{remark}

\begin{remark}\label{remequivSR}
It is a standard fact that the minimal time control problem for the
system (\ref{systeme}) with the constraint (\ref{contrainte}), is
equivalent to the sub-Riemannian problem associated to the $m$-tuple
of vector fields $(f_1,\ldots,f_m)$ (see \cite{Be} for a general
definition of a sub-Riemannian distance). In this context, there holds
$T_{\bar{x}}(x)=d_{SR}(\bar{x},x)$,
where $d_{SR}$ is the sub-Riemannian distance.
This implies that the functions $T_{\bar{x}}(\cdot)$ and
$d_{SR}(\bar{x},\cdot)$ share the same regularity properties.
In particular, the function $T_{\bar{x}}(\cdot)$ is continuous.
\end{remark}


\subsection{Class of controllers and notion of hybrid solution}
\label{contr_sect}
Let $f:\R^n\times \R^m\rightarrow \R^n$
be defined by $f(x,u)= \sum_{i=1}^m u_i f_i (x) .$
The system (\ref{systeme}) writes
\begin{equation}
    \label{ch:e1}
    \dot x(t) =f(x(t),u(t)) .
\end{equation}
Let $\bar{x}\in\R^n$ be fixed.

The controllers under consideration
in this paper depend on the continuous state $x\in\R^n$ and also on a
discrete variable $s_d\in \NNN$, where $\NNN$ is a nonempty
subset of $\N$.
According to the concept of a hybrid system of \cite{go-te}, we introduce the
following definition.

\begin{definition}
A hybrid feedback is a 4-tuple $(C,D,k,k_d)$, where
\begin{itemize}
\item $C$ and $D$ are subsets of $\R^n\times \NNN$;
\item $k:\R^n \times \NNN \rightarrow \R^m$ is a function;
\item $k_d: \R^n \times \NNN \rightarrow \NNN$ is a
function.
\end{itemize}
The sets $C$ and $D$ are respectively called the {\em controlled continuous
evolution set} and the {\em controlled discrete evolution set}.
\end{definition}

We next recall the notion of robustness to small noise (see
\cite{4}). Consider two functions $e$ and $d$ satisfying the
following {\em regularity assumptions}:
\begin{equation}\label{standing}
\begin{split}
& e(\cdot,\cdot), d(\cdot,\cdot)
 \in L^\infty_{loc}(\R^n\times [0,+\infty);\R^n),\\
&  e(\cdot,t), d(\cdot,t) \in C^0(\R^n,\R^n),\
\forall t\in[0,+\infty).
\end{split}
\end{equation}
We introduce these functions as a measurement noise~$e$
and an external disturbance~$d$.

Formally, the $k$-component of a hybrid feedback $(k,k_d,C,D)$ governs the
differential equation 
$$
\dot x  = f(x,k(x+e))+ d\ , \; \forall (x,s_d)\in C,
$$
whereas the $k_d$-component governs the jump equation
$$
s_d^+= k_d(x,s_d)\ , \; \forall (x,s_d)\in D .
$$
The set $C$ indicates where in the state space flow may occur
while the set
$D$ indicates where in the state space jumps may occur. The collection of
this flow equation and of this jump equation, under the perturbations $e$ and
$d$, is a perturbed hybrid system $\mathcal{H}_{(e,d)}$, as considered
e.g.\ in
\cite{nolcos04}. We next provide a precise definition of the notion of
solutions considered here.

This concept is
well studied in the literature (see \textit{e.g.}\
\cite{bm,bra,lnp,pat,chained,tav1}). Here, we consider the notion of solution
given in \cite{go-te,nolcos04}.
\begin{definition} \label{ch:eq:49}
Let $S=\bigcup_{j=0}^{J-1} [t_j,t_{j+1}]\times \{j\}$, where $J\in
\N\cup\{+\infty\}$ and $(x_0,s_0)\in \R^n\times
\NNN$. The domain $S$ is said to be a hybrid time domain.
A map $(x,s_d):
S \rightarrow \R^n\times \NNN$ is said to be a solution
of $\mathcal{H}_{(e,d)}$ with the initial condition $(x_0,s_0)$ if
\begin{itemize}
\item the map $x$ is continuous on $S$;
\item for every $j$, $0\leq j\leq J-1$, the map $x:t\in
(t_j,t_{j+1})\mapsto x(t,j)$ is absolutely continuous;
\item for every $j$, $0\leq j\leq J-1$ and
almost every $t\geq 0$, $(t,j)\in S$,
we have
\begin{equation}
\label{e135}
(x(t,j)+e(x(t,j),t),s_d(t,j))\in C,
\end{equation}
and
\begin{eqnarray}
\label{e134}
\dot x (t,j) &= &f(x(t),k(x(t,j)+e(x(t,j),t),s_d(t,j)))+d(x(t,j),t),
\\
\dot s_d(t,j)&=& 0 ;
\end{eqnarray}
(where the dot stands for the derivative with respect to the time
variable $t$)
\item for every $(t,j)\in S$, $(t,j+1)\in S$, we have
\begin{equation}
\label{e136}
(x(t,j)+e(x(t,j),t),s_d(t,j))\in D,
\end{equation}
and
\begin{eqnarray}
\label{e137}
x(t,j+1)&=&x(t,j),\\
s_d(t,j+1)&=& k_d(x(t,j)+e(x(t,j),t),s_d(t,j));
\end{eqnarray}
\item
$
(x(0,0),s_d(0,0))=(x_0,s_0)
$.
\end{itemize}
\end{definition}

In this context, we next define the concept of stabilization of
(\ref{ch:e1}) by a minimal time hybrid feedback law sharing a
robustness property with respect to measurement noise and external
disturbances. The usual Euclidean norm in ${\R}^n$ is denoted by
$|\cdot|$, and the open 
ball
centered at $\bar x$ with radius $R$ is denoted $B(\bar x,R)$. 
Recall that a function of class
${\mathcal{K}}_{\infty}$ is a function $\delta$:
$[0,+\infty)\rightarrow [0,+\infty)$ which is continuous,
increasing, satisfying $\delta(0)=0$ and $\lim_{R \rightarrow
+\infty}\delta(R)=+\infty$.

As usual, the system is said complete if all solutions are
maximally defined in $[0,+\infty)$ (see e.g. \cite{ABB}). More
precisely, we have the following definition.

\begin{definition}
Let $\rho:{\R}^n\rightarrow {\R}$ be a continuous function satisfying
\begin{equation}
\label{eq12}
\rho(x)>0 , \ \forall x\neq \bar{x} .
\end{equation}
We say that the {\em completeness assumption for $\rho$}
holds if,
for all $(e,d)$ satisfying the regularity assumptions (\ref{standing}),
and so that,
\begin{equation}\label{ch:eq31}
\sup_{[0,+\infty)}|e(x,\cdot)|\leq
\rho(x) , \
{\mathrm{esssup}}_{[0,+\infty)}|d(x,\cdot)|\leq\rho(x) , \ \forall x\in
{\R}^{n},
\end{equation}
for every $(x_0,s_{0})\in {\R}^{n}\times \NNN$,
there exists a maximal solution on $[0,+\infty)$
of $\mathcal{H}_{(e,d)}$ starting from $(x_0,s_0)$.
\end{definition}

Roughly speaking, the finite time convergence property means that
all solutions reach $\bar x$ within finite time. A precise definition of
this concept follows.

\begin{definition}
We say that the {\em uniform finite time convergence property} holds
if there exists a continuous function $\rho:{\R}^n\rightarrow {\R}$ satisfying
(\ref{eq12}), such that the completeness assumption
for $\rho$ holds, and if there exists
a function $\delta: [0,+\infty)\rightarrow [0,+\infty)$ of class
${\mathcal{K}}_\infty$ such that, for every $R>0$, there exists
$\tau=\tau(R)>0$, for all functions $e,d$
satisfying the
regularity assumptions (\ref{standing}) and inequalities (\ref{ch:eq31})
for this function $\rho$, for every
$x_0\in B(\bar{x},R)$, and every $s_0\in\NNN$, the
maximal solution $(x,s_d)$ of $\mathcal{H}_{(e,d)}$ starting from
$(x_0,s_0)$ satisfies
\begin{equation} \label{eq10}
|x(t,j)-\bar{x}|\leq \delta(R), \ \forall t\geq 0, \ (t,j)\in S,
\end{equation}
and
\begin{equation}
\label{eq11}
x(t,j)=\bar{x} , \ \forall t\geq \tau , \ (t,j)\in S.
\end{equation}
\end{definition}

We are now in position to introduce our main definition. It deals with 
closed-loop systems whose trajectories converge to the equilibrium within
quasi-minimal time and with a robustness property with respect to
measurement noise and external disturbances.

\begin{definition}\label{ch:def2}
The point $\bar{x}$ is said to be a
{\em semi-globally quasi-minimal time
robustly stabilizable equilibrium} for
the system (\ref{ch:e1}) if, for every
$\varepsilon>0$ and every compact subset $K\subset \R^n$, there exists
a hybrid feedback law $(C,D, k,k_d)$
satisfying the constraint
\begin{equation}
\label{eq3}
\Vert k(x,s_d)\Vert  \leq 1,
\end{equation}
where $\Vert\cdot\Vert$ stands for the Euclidean norm in ${\R}^m$,
such that:
\begin{itemize}
\item the uniform finite time convergence property
holds;
\item
there exists a continuous function $\rho_{\varepsilon,K}:{\R}^n\rightarrow {\R}$ satisfying
(\ref{eq12}) for $\rho=\rho_{\varepsilon,K}$, such that, for every
$(x_0,s_0)\in K\times \NNN$, all functions $e,d$ satisfying the
regularity assumptions (\ref{standing}) and inequalities (\ref{ch:eq31})
for $\rho=\rho_{\varepsilon,K}$,
the maximal solution of $\mathcal{H}_{(e,d)}$ starting from
$(x_0,s_0)$ reaches $\bar{x}$ within time
$T_{\bar{x}}(x_0)+\varepsilon $, where $T_{\bar{x}}(x_0)$ denotes
the minimal time to steer the system (\ref{ch:e1}) from $x_0$ to $\bar{x}$,
under the constraint $\Vert u\Vert \leq 1$.
\end{itemize}
\end{definition}


\subsection{Main result}\label{ch:main_sect}
The main result of this article is the following.
\begin{theorem}\label{CH:TH1}
Let $\bar{x}\in \R^n$.
If there exists no nontrivial minimal time singular trajectory of
(\ref{systeme}), (\ref{contrainte}), starting from $\bar{x}$,
then $\bar{x}$ is a semi-globally quasi-minimal time robustly
stabilizable equilibrium for the system (\ref{systeme}), under the
constraint (\ref{contrainte}).
\end{theorem}

\begin{remark}
The problem of global quasi-minimal time robust stabilization
(\textit{i.e.}\ $K=\R^n$ in Definition \ref{ch:def2}) cannot be
achieved a priori because measurement noise may then accumulate
and slow down the solution reaching $\bar{x}$ (compare with
\cite{Bressantorino}).
\end{remark}

\begin{remark}
The assumption of the absence of nontrivial singular minimizing
trajectory is crucial.
Notice the following facts, which show the relevance of this assumption:
\begin{itemize}
\item
if $m\geq n$ and if the vector fields $f_1,\ldots,f_m$, are everywhere
linearly independent, then there exists no
singular trajectory. In this case, we are actually in the framework
of Riemannian geometry (see Remark \ref{remequivSR}).
\item
Let ${\mathcal F}_m$ be the set of $m$-tuples of linearly
independent vector
fields
$(f_1,\ldots,f_m)$, endowed with the $C^\infty$ Whitney
topology. If $m\geq 3$, there exists an open dense subset of ${\mathcal
F}_m$, such that any control system of the form (\ref{systeme}),
associated to a $m$-tuple of this subset, admits no nontrivial
singular minimizing trajectory (see \cite{CJTcras,CJT}, see also
\cite{AG} for the existence of a dense set only).
Hence generically the conclusion of Theorem \ref{CH:TH1} holds
without assuming the absence of nontrivial singular minimizing
trajectories.
\item
If there exist singular minimizing trajectories, then the conclusion
on subanalyticity of the function $T$ may fail (see \cite{BC,these}),
and we cannot a priori prove that the set $\mathcal{S}$ of
singularities of $T$ is a stratifiable manifold, which is the crucial
fact in order to define a hybrid strategy.
\end{itemize}
\end{remark}


\subsection{Short description of the proof}\label{secshortdescription}
The strategy of the proof of Theorem \ref{CH:TH1} is the following.

\begin{figure}[ht]
\centerline{\input{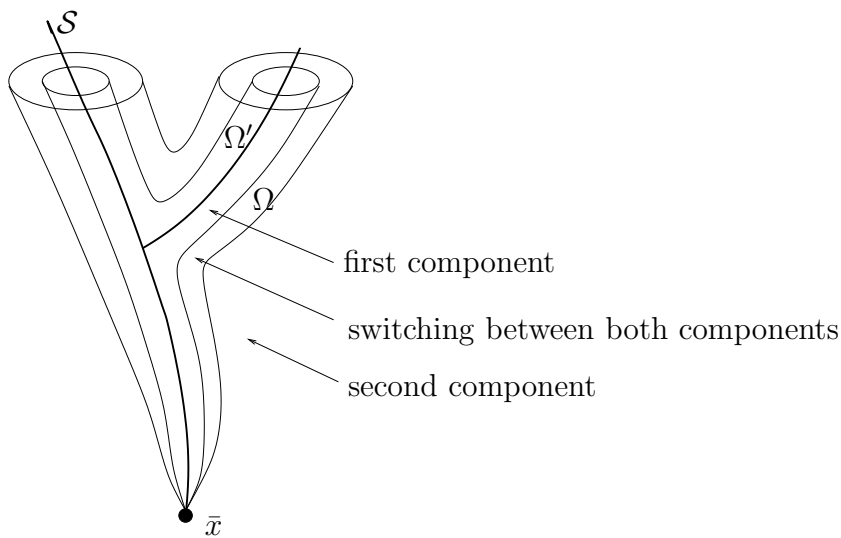}} \caption{Switching
strategy.}\label{fighybrid}
\end{figure}

Under the assumption of the absence of nontrivial singular minimal time
trajectory, the minimal time function $T_{\bar{x}}$ associated to the
system (\ref{systeme}), (\ref{contrainte}), is subanalytic, and hence,
is analytic outside a stratified submanifold $\mathcal{S}$ of $\R^n$,
of codimension greater than or equal to one. Therefore,
the corresponding minimal time feedback controller (further precisely
defined in Section \ref{secdefopt}) is continuous (even
analytic) on $\R^n\setminus{\mathcal{S}}$ (see Figure
\ref{fighybrid}).
In a neighborhood of $\mathcal{S}$, it is therefore
necessary to use other controllers, and then to define an adequate
switching strategy.

More precisely, the proof of Theorem \ref{CH:TH1} relies on both
following key lemmas.
\begin{lemma}\label{e160}
For every
$\eps>0$, there exists a neighborhood $\Omega$ of $\mathcal{S}$ such that,
for every stratum%
\footnote{Since $\mathcal{S}$ is a stratified submanifold of $\R^n$ of
codimension greater than or equal to one, there exists a partition
$(M_i)_{i\in \N}$ of $\mathcal{S}$, where $M_i$ is a stratum,
\textit{i.e.}, a locally closed submanifold of $\R^n$. Recall that, if
$M_i\cap \partial M_j \neq \emptyset$, then $M_i\subset M_j$ and
$\mathrm{dim}(M_i) < \mathrm{dim}(M_j)$.} %
 $M_i$ of $\mathcal{S}$,
there exist a nonempty subset $\NNN_i$ of $\N$, a locally
finite family $(\Omega_{i,p})_{p \in \NNN_i}$ of open subsets of $\Omega$,
a sequence of smooth controllers
$u_{i,p}$ defined in a neighborhood of $\Om_{i,p}$,
satisfying $\Vert u_{i,p}\Vert \leq 1$, and there exists
a continuous function
$\rho_{i,p}:\R^n\ra [0,+\infty)$ satisfying $\rho_{i,p}(x)>0$ whenever $x\neq
\bar{x}$, such that every solution of
\begin{equation}\label{e163}
\dot x(t)= f(x(t),u_{i,p}(x(t)+e(x(t),t)))+d(x(t),t),
\end{equation}
where $e,d: \R^n \times [0,+\infty) \rightarrow \R^n$ are two functions
satisfying the regularity assumptions (\ref{standing}) and
\begin{equation}
 \sup_{[0,+\infty)}|e(x,\cdot)|\leq \rho_{i,p}(x)
,\ \
{\mathrm{esssup}}_{[0,+\infty)}|d(x,\cdot)|\leq \rho_{i,p} (x) ,
\end{equation}
starting from $\Om_{i,p}$ and maximally defined on $[0,T)$, leaves
$\Om$ within time $\eps$; moreover, there exists a
function $\delta_{i,p}$ of class ${\mathcal{K}}_{\infty}$ such that,
for every $R>0$, every such solution starting from $\Om_{i,p}\cap
B(\bar{x},R)$ satisfies 
\begin{equation}\label{deltaip}
|x(t)-\bar{x}|\leq  \delta_{i,p}(R) ,\ \  \fa t\in [0,T) .
\end{equation}
\end{lemma}

According to this lemma, in a neighborhood $\Omega$ of $\mathcal{S}$,
there exist controllers steering the system outside
$\Omega$. Moreover, since this neighborhood can be
chosen arbitrarily thin, the time $\varepsilon$ needed for its
traversing is arbitrarily small.

Outside $\Omega$, the optimal controller is analytic. The following
lemma shows that this controller shares an invariance property; in
brief, it gives rise to trajectories never crossing again the singular
set $\mathcal{S}$.

\begin{lemma}\label{leminvariance2}
For every neighborhood $\Omega$ of $\mathcal{S}\setminus\{ \bar{x}
\}$ in $\R^n$, there exists a neighborhood $\Omega'$ of
$\mathcal{S}\setminus\{ \bar{x} \}$ in $\R^n$, satisfying
\begin{equation}\label{17juin}
\Omega'\subsetneq\clo (\Omega')\subsetneq \Omega,
\end{equation}
such that every trajectory of the closed-loop system (\ref{systeme})
with the optimal controller, starting from any point $x\in
\R^n\setminus\Omega$, reaches $\bar{x}$ in minimal time, and is
contained in $\R^n\setminus\Omega'$.
\end{lemma}

\medskip

Finally, our hybrid strategy is the following. For every $\eps>0$,
there exists a neighborhood $\Omega$ of the singular set
$\mathcal{S}$, and there exist controllers which steer the system outside
this neighborhood in time less than $\eps$. Outside $\Omega$, there
exists a continuous controller and giving rise to trajectories never
crossing again $\mathcal{S}$ and joining $\bar{x}$ in minimal time.

It is therefore necessary to define an adequate switching strategy
connecting both controllers (see Figure \ref{fighybrid}).
This is achieved in the context of hybrid
systems, using an hysteresis strategy. The
first component consists of controllers which are defined in $\Omega$, and
whose existence is stated in
Lemma \ref{e160}. The second component of the
hysteresis is defined by the optimal controller, outside $\Omega$; 
Both components are united using an hysteresis, by
adding a dynamical discrete variable $s_d$ and using a hybrid feedback
law. With this resulting hybrid controller, the time needed to join
$\bar{x}$, from any point $x_0$ of $\R^n$, is less than
$T_{\bar{x}}(x_0)+\varepsilon$.

\medskip

The next section, devoted to the detailed proof of Theorem
\ref{CH:TH1}, is organized as follows.

The first component of the hysteresis is defined in Section
\ref{secsecondcomponent}, and Lemma \ref{e160} is proved.

Section \ref{firstcomponent} concerns the definition and properties of
the second component of the hysteresis, defined by the minimal time
controller. In Section \ref{seccomputation}, we recall how to
compute minimal time trajectories of the system (\ref{systeme}),
(\ref{contrainte}), using the Pontryagin Maximum Principle. We then
provide in Section \ref{seccutlocus} a crucial result on the cut locus
(Proposition \ref{propcutlocus}). The optimal feedback controller is
defined in Section \ref{secdefopt}; basic facts on subanalytic
functions are recalled, permitting to define the singular set
$\mathcal{S}$. Invariance properties of this optimal controller are
then investigated: Lemma \ref{leminvariance2} is proved in Section
\ref{secproofleminvariance2}; robustness properties are given and
proved in Section \ref{secoptrobust}.

The hybrid controller is then defined in Section \ref{sechybrid}.
A definition of a hybrid control system, and properties of solutions,
are given in Sections \ref{secdefgeneralhybrid} and
\ref{secpropertiessolutionshybrid}.
A precise description of the switching strategy is provided in Section
\ref{secswitching}.
Theorem \ref{CH:TH1} is proved in Section \ref{e157}.


\section{Proof of Theorem \ref{CH:TH1}}\label{secproofmainthm}
In what follows, let $\bar{x}\in \R^n$  be fixed.


\subsection{The first component of the hysteresis}
\label{secsecondcomponent}
The first component of the hysteresis consists of a set of controllers,
defined in a neighborhood of $\mathcal{S}$, whose
existence is stated in Lemma \ref{e160}.
Hereafter, we provide a proof of this lemma.

\begin{proof}[Proof of Lemma \ref{e160}]
First of all, recall that, on the one hand, the minimal time function
coincides with the sub-Riemannian distance associated to the $m$-tuple
$(f_1,\ldots,f_m)$ (see Remark \ref{remequivSR}); on the other,
since the Lie Algebra Rank Condition holds,
the topology defined by the sub-Riemannian distance $d_{SR}$ coincides
with the Euclidean topology of $\R^n$, and, since
$\R^n$ is complete, any two points of $\R^n$ can be joined
by a minimizing path (see \cite{Be}).

Let $\eps > 0$ fixed. Since $\mathcal{S}$ is a stratified submanifold
of $\R^n$, there exists a neighborhood $\Omega$ of $\mathcal{S}$
satisfying the following property: for every $y\in\mathcal{S}$, there
exists $z\in\R^n\setminus\clo(\Omega)$ such that
$d_{SR}(y,z)<\eps$.

Consider a stratum $M_i$ of $\mathcal{S}$.
For every $y\in M_i$, let $z\in\R^n\setminus\clo(\Omega)$ such that
$d_{SR}(y,z)<\eps$.
The Lie Algebra Rank Condition implies that there exists an
open-loop control $t\mapsto u_y(t)$, defined on $[0,T)$ for
a $T>\eps$, satisfying the constraint $\Vert u_y \Vert \leq 1$, such
that the associated trajectory $x_y(\cdot)$ (which can be assumed to
be one-to-one), solution of (\ref{systeme}),
starting from $y$, reaches $z$ (and thus, leaves
$\clo(\Omega)$) within time $\eps$.
Using a density argument, the control $u_y$ can be moreover
chosen as a smooth function (see \cite[Theorem 2.8 p.~21]{Be} for the
proof of this statement).
Since the trajectory is one-to-one, the
open-loop control $t\mapsto u_y(t)$ can be considered as a feedback
$t\mapsto u_y (x_y (t))$ along $x_y(\cdot)$.
Consider a smooth extension of $u_y$ on $\R^n$,
still denoted $u_y$, satisfying the constraint $\Vert
u_{y}(x)\Vert \leq 1$, for every $x\in\R^n$. By continuity, there
exists a neighborhood $\Om_y$ of $y$, and positive real numbers
$\delta_y$ and $\rho_y$, such that every solution of
\begin{equation}\label{e163y}
\dot x(t)= f(x(t),u_{y}(x(t)+e(x(t),t)))+d(x(t),t),
\end{equation}
where $e,d: \R^n \times [0,+\infty) \rightarrow \R^n$ are two functions
satisfying the regularity assumptions (\ref{standing}) and
\begin{equation*}
 \sup_{[0,+\infty)}|e(x,\cdot)|\leq \rho_y
,\ \
{\mathrm{esssup}}_{[0,+\infty)}|d(x,\cdot)|\leq \rho_y \ ,
\end{equation*}
starting from $\Om_y$ and maximally defined on $[0,T)$, leaves
$\Om$ within time $\eps$; moreover,
$$
|x(t)-\bar{x}|\leq  \delta_y ,\  \fa t\in [0,T) .
$$
Repeat this construction for each $y\in M_i$.

Now, on the one hand, let $(y_p)_{p\in\NNN_i}$ be a sequence of points of
$M_i$ such that
the family $(\Om_{y_p})_{p\in\NNN_i}$ is a locally finite
covering of $M_i$, where $\NNN_i$ is a subset of $\N$.
Define $\Om_{i,p}=\Om_{y_p}$ and $u_{i,p}=u_{y_p}$.

On the other hand, the existence of a continuous function
$\rho_{i,p}:\R^n\rightarrow [0,+\infty)$, satisfying $\rho_{i,p}(x)>0$
whenever $x\neq x$, follows for the continuity of
solutions with respect to disturbances. 
The existence of a function $\delta_{i,p}$ of class
${\mathcal{K}}_{\infty}$ such that (\ref{deltaip}) holds is obvious.

Repeat this construction for every stratum $M_i$ of $\mathcal{S}$.
Then, the statement of the lemma follows.
\end{proof}
\subsection{The second component of the hysteresis}\label{firstcomponent}
\subsubsection{Computation of minimal time trajectories}\label{seccomputation}
Let $x_1\in\R^n$, and $x(\cdot)$ be a minimal time
trajectory, associated to a control $u(\cdot)$, steering the system
(\ref{systeme}), (\ref{contrainte}), from $\bar{x}$ to $x_1$, in time
$T=T_{\bar{x}}(x_1)$. According to Pontryagin's maximum principle (see
\cite{P}), the trajectory $x(\cdot)$ is the projection of an
\textit{extremal}, \textit{i.e.}, a triple
$(x(\cdot),p(\cdot),u(\cdot))$ solution of the constrained Hamiltonian
system
$$\dot{x}(t)=\frac{\partial H}{\partial p}(x(t),p(t),u(t)),\ \
\dot{p}(t)=-\frac{\partial H}{\partial x}(x(t),p(t),u(t)),$$
$$H(x(t),p(t),p^0,u(t))=\max_{\Vert v\Vert\leq 1}H(x(t),p(t),p^0,v),$$
almost everywhere on $[0,T]$, where
$$
H(x,p,u)= \sum_{i=1}^m u_i\langle p, f_i(x)\rangle
$$
is the \textit{Hamiltonian} of the optimal control problem, and
$p(\cdot)$ (called \textit{adjoint vector}) is an absolutely continuous
mapping on $[0,T]$ such that $p(t)\in \R^n\setminus\{0\}$.
Moreover, the function $t\mapsto \max_{\Vert v\Vert\leq
  1}H(x(t),p(t),p^0,v)$ is Lipschitzian, and
everywhere constant on $[0,T]$. If this constant is not equal to zero,
then the extremal is said \textit{normal}; otherwise it is said
\textit{abnormal}.

\begin{remark}\label{rem1.4}
Any singular trajectory is the projection of an abnormal extremal,
and conversely.
\end{remark}

Controls associated to normal extremals can be computed as
\begin{equation}\label{dennonnul}
u_i(t)=\frac{\langle p(t),f_i(x(t)\rangle}
{\sqrt{{\sum_{j=1}^m\langle p(t),f_j(x(t)\rangle^2}}},\ i=1,\ldots,m.
\end{equation}
Indeed, note that, by definition of normal extremals, the denominator
of (\ref{dennonnul}) cannot vanish. It follows that
normal extremals are solutions of the Hamiltonian system
\begin{equation}\label{reduit}
\dot{x}(t)=\frac{\partial H_1}{\partial p}(x(t),p(t)),\
\
\dot{p}(t)=-\frac{\partial H_1}{\partial x}(x(t),p(t)),
\end{equation}
where
$$
H_1(x,p) = \sqrt{\sum_{i=1}^m \langle p,f_i(x)\rangle^2}.
$$
Notice that $H_1(x(t),p(t))$ is constant, nonzero, along each normal
extremal. Since $p(0)$ is defined up to a multiplicative scalar, it is
usual to normalize it so that $H_1(x(t),p(t))=1$.
Hence, we introduce the set
$$ X = \{ p_0\in \R^n\ \vert\ H_1(\bar{x},p_0)=0 \}.$$
It is a submanifold of $\R^n$ of codimension one, since
$\frac{\partial H_1}{\partial p}(\bar{x},p_0) = \dot{x}(0)\neq
0$ (see \cite{BCT} for a similar construction in the general case).
There exists a connected open subset $U$ of $[0,+\infty)\times X$ such that,
for every $(t^*,p_0)\in X$, the differential system
(\ref{reduit}) has a well defined smooth solution on $[0,t^*]$ such that
$x(0)=\bar{x}$ and $p(0)=p_0$.

\begin{definition}\label{defexp}
The smooth mapping
$$\begin{array}{rrcl}
\textrm{exp}_{\bar{x}}:&       U & \longrightarrow & \R^n    \\
                     &   (t,p_0)   & \longmapsto     & x(t)
\end{array}$$
where $(x(\cdot),p(\cdot))$ is the solution of the system (\ref{reduit})
such that $x(0)=\bar{x}$ and $p(0)=p_0\in X$, is called
\textit{exponential mapping} at the point $\bar{x}$.
\end{definition}

The exponential mapping parameterizes normal extremals.
Note that the domain of $\textrm{exp}_{\bar{x}}$ is a subset of
$\R\times X$ which is locally diffeomorphic to $\R^n$ (since we are in
the normal case).

\begin{definition}
A point $x\in\textrm{exp}_{\bar{x}}(U)$ is said to be \textit{conjugate}
to $\bar{x}$ if it is a critical value of the mapping $\textrm{exp}_{\bar{x}}$,
\textit{i.e.}, if there exists $(t_c,p_0)\in U$ such that
$x=\textrm{exp}_{\bar{x}}(t_c,p_0)$ and the differential
$d\,\textrm{exp}_{\bar{x}}(t_c,p_0)$ is not onto. The \textit{conjugate
locus} of $\bar{x}$, denoted by ${\mathcal{C}}(\bar{x})$, is defined as the
set of all points conjugate to $\bar{x}$.

With the previous notations, define ${\mathcal{C}}_{min}(\bar{x})$ as the
set of points $x\in {\mathcal{C}}(\bar{x})$ such that the trajectory
$t\mapsto \textrm{exp}_{\bar{x}}(t,p_0)$ is minimizing between $\bar{x}$ and
$x$.
\end{definition}


\subsubsection{The cut locus}\label{seccutlocus}
A standard definition is the following.
\begin{definition}
A point $x\in \R^n$ is not a \textit{cut point} with respect to
$\bar{x}$ if there exists a minimal time trajectory of (\ref{systeme}),
(\ref{contrainte}), joining $\bar{x}$ to $x$,
which is the strict restriction of a minimal time trajectory starting
from $\bar{x}$. The \textit{cut locus} of $\bar{x}$, denoted by
${\mathcal{L}}(\bar{x})$, is defined as the set of all cut points with
respect to $\bar{x}$.
\end{definition}

In other words, a cut point is a point at which a minimal time
trajectory ceases to be optimal.

\begin{remark}\label{remconjanalytic}
In the analytic case, it follows from the theory of conjugate
points that every nonsingular minimal time trajectory ceases to be
minimizing beyond its first conjugate point (see for instance
\cite{AgSach,BC}). Hence, if there exists no
singular minimal time trajectory starting from $\bar{x}$, then
${\mathcal{C}}_{min}(\bar{x})\subset {\mathcal{L}}(\bar{x})$.
\end{remark}

The following result on the cut locus is crucial for the proof of Theorem
\ref{CH:TH1}.

\begin{proposition}\label{propcutlocus}
Assume that the vector fields $f_1,\ldots,f_m$ are analytic,
and that there exists no singular minimal time trajectory starting
from $\bar{x}$. Then the set of points of $\R^n$ where the
function $T_{\bar{x}}(\cdot)$ is not analytic is equal to the
cut locus of $\bar{x}$, that is
\begin{equation}
\mathrm{Sing}\ T_{\bar{x}}(\cdot)={\mathcal{L}}(\bar{x}).
\end{equation}
\end{proposition}

\begin{remark}
Under the previous assumptions, one can prove that the set of
points of $\R^n$ where $T_{\bar{x}}(\cdot)$ is analytic is equal to
the set of points where $T_{\bar{x}}(\cdot)$ is of class $C^1$.
\end{remark}

\begin{proof}
Let $x\in \R^n$ so that $T_{\bar{x}}(\cdot)$ is analytic at $x$.
Then there exists a neighborhood $V$ of $x$ in $\R^n$ such that
$T_{\bar{x}}(\cdot)$ is analytic on $V$. Let us prove that $x\notin
{\mathcal{L}}(\bar{x})$. It follows from the maximum principle and the
Hamilton-Jacobi theory (see \cite{P}) that, for every $y\in V$,
there exists a unique minimal time trajectory joining $\bar{x}$ to $y$,
having moreover a normal extremal lift
$(x(\cdot),p(\cdot),u(\cdot))$ satisfying
$$p(T_{\bar{x}}(y)) = \nabla T_{\bar{x}}(y)$$
(compare with \cite[Proposition 2.3]{RiffordTrelat}).
Set $ U_1 = \textrm{exp}_{\bar{x}}^{-1}(V)$. It follows easily from
Cauchy-Lipschitz Theorem that the mapping $\textrm{exp}_{\bar{x}}$ is an
analytic diffeomorphism from $U_1$ into $V$.
Hence, obviously, the point $x$ is not in the cut locus of $\bar{x}$.

\medskip

Conversely, let $x\notin {\mathcal{L}}(\bar{x})$. To prove that
$T_{\bar{x}}(\cdot)$ is analytic at $x$, we need the two following lemmas.

\begin{lemma}\label{lem2prop}
The point $x$ is not conjugate to $\bar{x}$, and is joined from
$\bar{x}$ by a unique minimal time trajectory.
\end{lemma}

\begin{proof}[Proof of Lemma \ref{lem2prop}.]
 From the assumption of the absence
of singular minimal time trajectory, there exists a nonsingular
minimal time trajectory joining $\bar{x}$ to $x$. From Remark
\ref{remconjanalytic}, the point $x$ is not conjugate to $\bar{x}$.

By contradiction, suppose that $x$ is joined from $\bar{x}$ by at least
two minimal time trajectories. By assumption, these two
trajectories must admit normal extremal lifts. Since the structure is
analytic, their junction at the point $x$ is necessarily not
smooth. This implies that both trajectories loose their optimality at the
point $x$ (indeed if not, there would exist a nonsmooth normal
extremal, which is absurd), and thus $x\in
{\mathcal{L}}(\bar{x})$. This is a contradiction.
\end{proof}

\begin{lemma}\label{lem3prop}
There exists a neighborhood $V$ of $x$ in $\R^n$ such that every
point $y\in V$ is not conjugate to $\bar{x}$, and there exists a
unique (nonsingular) minimal time trajectory joining $\bar{x}$ to $y$.
\end{lemma}

\begin{proof}[Proof of Lemma \ref{lem3prop}]
Let $p_0\in X$ so that $x=\textrm{exp}_{\bar{x}}(T_{\bar{x}}(x),p_0)$. Since $x$
is not conjugate to $\bar{x}$, the exponential mapping
$\textrm{exp}_{\bar{x}}$ is a diffeomorphism from a neighborhood $U_1$
of $(T_{\bar{x}}(x),p_0)$ in $U$ into a neighborhood $V$ of $x$ in $\R^n$.
Set $U_2=\textrm{exp}_{\bar{x}}^{-1}(V)$.

Let us prove that $\textrm{exp}_{\bar{x}}$
is proper from $U_2$ into $V$.
We argue by contradiction, and suppose that there exists a sequence
$(x_n)_{n\in\N}$ of points of $V$ converging towards $x$, such that
for each integer $n$ there exists $p_n\in X$, satisfying
$(T_{\bar{x}}(x_n),p_n)\in U_2$ and
$x_n=\textrm{exp}_{\bar{x}}(T_{\bar{x}}(x_n),p_n)$, such that $(p_n)_{n\in\N}$ is not
bounded. It then follows from \cite[Lemmas 4.8 and 4.9]{trelatJDCS}
(see also \cite[Fact 1 p.~378]{trelatIHP}) that $x$ is joined from
$\bar{x}$ by a singular control $u$. In particular, $x$ is conjugate to
$\bar{x}$; this is a contradiction.

Therefore, the set $\{p\ \vert\ \textrm{exp}_{\bar{x}}(T_{\bar{x}}(x),p)=x\}$
is compact in $U_2$. Moreover, since $x$ is not conjugate to $\bar{x}$, this
set has no cluster point, and thus is finite. As a consequence, up to
reducing $V$, we assume that
$V$ is a connected open subset of $\textrm{exp}_{\bar{x}}(U_2)$, and that
$U_2$ is a finite union of disjoint connected open sets, all of which being
diffeomorphic to $V$ by the mapping $\textrm{exp}_{\bar{x}}$. We infer
that every point $y \in V$ is not conjugate to $\bar{x}$.
Hence, the mapping $\textrm{exp}_{\bar{x}}$ is a proper submersion from
$U_2$ into $V$, and thus is a fibration with finite degree. Since,
from Lemma \ref{lem2prop}, there exists
a unique minimal time trajectory joining $\bar{x}$ to $x$, this degree is
equal to one, that is, $\textrm{exp}_{\bar{x}}$ is a diffeomorphism from
$U_2$ into $V$. The conclusion follows.
\end{proof}

It follows from the previous lemma that
$ (T_{\bar{x}}(y),p_0) =  \textrm{exp}_{\bar{x}}^{-1}(y)$,
for every $y\in V$, and hence $T_{\bar{x}}(\cdot)$ is analytic on $V$.
\end{proof}


\subsubsection{Definition of the optimal controller}\label{secdefopt}
By assumption, there does not exist any nontrivial
singular minimal time trajectory starting from $\bar{x}$.
Under these assumptions, the function $T_{\bar{x}}(\cdot)$ is
\textit{subanalytic} outside $\bar{x}$ (see \cite{Ag,AG,these},
combined with Remark \ref{remequivSR}).

\medskip

For the sake of completeness, we recall below the definition of a
subanalytic function (see \cite{Hardt,Hironaka}), and some
properties that are used in a crucial way in the present paper (see
\cite{Tamm}).

Let $M$ be a real analytic finite dimensional manifold. A subset $A$
of $M$ is said to be {\it{semi-analytic}} if and only if, for every
$x\in M$, there exists a neighborhood $U$ of $x$ in $M$ and $2pq$
analytic functions $g_{ij}, h_{ij}$ ($1\leq i\leq p$ and $1\leq j\leq
q$), such that
$$A\cap U=\bigcup_{i=1}^p\{y\in U\ \vert\ g_{ij}(y)=0\ \mathrm{and}\
h_{ij}(y)>0, \ j=1\ldots q\} .  $$
Let SEM($M$) denote the set of semi-analytic subsets of $M$.
The image of a semi-analytic subset by a proper analytic mapping is
not in general semi-analytic, and thus this class has to be enlarged.

A subset $A$ of $M$ is said to be {\it{subanalytic}} if and only if,
for every $x\in M$, there exist a neighborhood $U$ of $x$ in $M$ and
$2p$ couples $(\Phi_i^\delta,A_i^\delta)$ ($1\leq i\leq p$ and
$\delta=1,2$), where $A_i^\delta\in SEM(M_i^\delta)$, and where
the mappings $\Phi_i^\delta~: M_i^\delta \rightarrow M$ are proper
analytic, for real analytic manifolds $M_i^\delta$,  such that
$$A\cap U=\bigcup_{i=1}^p(\Phi_i^1(A_i^1)\backslash \Phi_i^2(A_i^2)) .  $$
Let SUB($M$) denote the set of subanalytic subsets of $M$.

The subanalytic class is closed by union, intersection, complementary,
inverse image by an analytic mapping, image by a proper analytic
mapping. In brief, the subanalytic class is {\it{o-minimal}} (see
\cite{DriesMiller}). Moreover subanalytic sets are {\it{stratifiable}}
in the following sense. A {\it{stratum}} of a
differentiable manifold $M$ is a locally closed sub-manifold of $M$.
A locally finite partition ${\mathcal S}$ of $M$ is a
{\it{stratification}} of $M$ if any $S\in{\mathcal S}$ is a stratum
such that
$$\forall T\in {\mathcal S}\quad T\cap \partial S\neq\emptyset\Rightarrow
T\subset \partial S\ \mathrm{and}\ \mathrm{dim}\ T<\mathrm{dim}\ S  . $$

Finally, a mapping $f:M\rightarrow N$ between two analytic manifolds
is said to be {\it{subanalytic}} if its graph is a subanalytic subset
of $M\times N$.

Let $M$ be an analytic manifold, and $f$ be a subanalytic function on
$M$. The \textit{analytic singular support} of $f$
is defined as the complement of the set of points $x$
in $M$ such that the restriction of $f$ to some neighborhood of $x$ is
analytic. The following property is of great interest in the present
paper (see \cite{Tamm}): the analytic singular support of $f$ is
subanalytic (and thus, in particular, is stratifiable). If
$f$ is moreover  locally bounded on $M$, then it is moreover of
codimension greater than or equal to one.

\medskip

Turn back to our problem. The function $T_{\bar{x}}(\cdot)$ is subanalytic
  outside $\bar{x}$, and hence, its singular set
${\mathcal{S}} = \mathrm{Sing}\ T_{\bar{x}}(\cdot)$ (\textit{i.e.},
the analytic singular support of $T_{\bar{x}}(\cdot)$)
is a stratified submanifold of $\R^n$, of codimension greater than
or equal to $1$.

\begin{remark}
Note that the point $\bar{x}$ belongs to the adherence of
${\mathcal{S}}$ (see \cite{Ag}).
\end{remark}

Outside the singular set $\mathcal{S}$, it follows from the dynamic
programming principle (see \cite{P}) that the minimal time
controllers steering a point $x\in \R^n\setminus{\mathcal{S}}$ to
$\bar{x}$ are given by the closed-loop formula
\begin{equation}\label{eq20}
u_i(x) = -\frac{\langle \nabla T_{\bar{x}}(x),f_i(x)\rangle}
{\sqrt{\sum_{j=1}^m \langle \nabla T_{\bar{x}}(x),f_j(x)\rangle^2}},\
i=1,\ldots,m.
\end{equation}

The objective is to construct neighborhoods of
${\mathcal{S}}\setminus\{ \bar{x}\}$ in $\R^n$ whose complements
share invariance properties for the optimal flow. This is the contents
of Lemma \ref{leminvariance2}, proved next.


\subsubsection{Proof of Lemma \ref{leminvariance2}}\label{secproofleminvariance2}
It suffices to prove that, for every compact subset $K$ of $\R^n$,
for every neighborhood $\Omega$ of $\mathcal{S}\setminus\{ \bar{x}
\}$ in $\R^n$, there exists a neighborhood $\Omega'$ of
$\mathcal{S}\setminus\{ \bar{x} \}$ in $\R^n$, satisfying (\ref{17juin}),
such that every trajectory of the
closed-loop system (\ref{systeme}) with the optimal controller,
joining a point $x\in (\R^n\setminus\Omega)\cap K$ to $\bar{x}$,
is contained in $\R^n\setminus\Omega'$.

By definition of the cut locus, and using Proposition
\ref{propcutlocus}, every optimal trajectory joining a point $x\in
(\R^n \setminus\Omega)\cap K$ to $\bar{x}$ does not intersect
$\mathcal{S}$, and thus has a positive distance to the set
$\mathcal{S}$. Using the assumption of the absence of nontrivial
singular minimizing trajectories starting from $\bar{x}$, a
reasoning similar to the proof of Lemma \ref{lem3prop} proves that
the optimal flow joining points of
the compact set $(\R^n\setminus\Omega)\cap K$ to $\bar{x}$ is
parameterized by a compact set. Hence, there exists a positive real
number $\delta>0$ so that every optimal trajectory joining a point
$x\in (\R^n\setminus\Omega)\cap K$ to $\bar{x}$ has a distance to
the set $\mathcal{S}$ which is greater than or equal to $\delta$.
The existence of $\Omega'$ follows.


\subsubsection{Robustness properties of the optimal controller}\label{secoptrobust}
In this section, we prove robustness properties of the Carath\'eodory
solutions of system (\ref{ch:e1}) in closed-loop with this feedback
optimal controller. Given $e,d: {\R}^n\times
[0,+\infty) \rightarrow{\R}^{n}$, the perturbed system in
closed-loop with the optimal controller (denoted $u_{opt}$) writes
\begin{equation}
\label{ch:eq:5} \dot x(t) =f(x(t),u_{opt}(x(t)+e(x(t),t)))
+d(x(t),t) .
\end{equation}
Since the optimal controller is continuous outside the singular
set $\mathcal{S}$, it enjoys a natural robustness property, stated
below. In the next result, the notation $d(x,\mathcal{S})$ stands for the
Euclidean distance from $x$ to $\mathcal{S}$.

\begin{lemma}\label{ch:le1}
There exist a continuous
function $\rho_{opt}:\R\rightarrow {\R}$ satisfying
\begin{equation}
\label{ch:eq:16bis}
\rho_{opt}(\xi)>0   , \ \forall \xi\neq 0  ,
\end{equation}
and a continuous function $\delta_{opt}:[0,+\infty)\rightarrow
[0,+\infty)$ of class ${\mathcal{K}}_{\infty}$ such that the
following three properties hold:
\begin{itemize}
\item {\em Robust Stability}

For every neighborhood $\Om$ of $\mathcal{S}$, there exists a
neighborhood $\Om'\subset \Om$ of $\mathcal{S}$, such that, for
all $e,d: \R^n\times [0,+\infty) \rightarrow \R^n$ satisfying the
regularity assumptions (\ref{standing}) and, for every $x\in\R^n$,
\begin{equation}
\label{ch:eq:20}
\sup_{[0,+\infty)}|e(x,\cdot)|\leq \rho_{opt}(d(x,\mathcal{S}))
,\ {\mathrm{esssup}}_{[0,+\infty)}|d(x,\cdot)|\leq \rho_{opt}(d(x,\mathcal{S}))
 ,
\end{equation}
and for every ${x}_0\in \R^n\setminus \Om$, there exists a unique
Carath\'eodory solution $x(\cdot)$ of (\ref{ch:eq:5}) starting
from ${x}_0$, maximally defined on $[0,+\infty)$, and satisfying
$x(t)\in\R^n\setminus \Om'$, for every $t> 0$.
\item {\em Finite time convergence}

For every $R>0$, there exists $\tau_{opt}=\tau_{opt}(R)>0$
such that, for all $e,d:
\R^n\times [0,+\infty) \rightarrow \R^n$ satisfying the regularity assumptions
(\ref{standing}) and (\ref{ch:eq:20}), for every
$x_0 \in \R^n$ with $|x_0-\bar{x}|\leq R$, and every
maximal solution $x(\cdot)$ of (\ref{ch:eq:5}) starting from
$x_0$, one has
\begin{equation} \label{ch:e27}
|x(t)-\bar{x}|\leq \delta_{opt}(R) ,\; \forall t\geq 0  ,
\end{equation}
\begin{equation}
\label{eq2}
x(t)=\bar{x} ,\; \forall t\geq \tau_{opt}  ,
\end{equation}
and
\begin{equation}
\label{eq6}
\|u_{opt}(x(t))\|\leq 1 ,\; \forall t\geq 0  .
\end{equation}
\item {\em Optimality}

For every neighborhood $\Om$ of $\mathcal{S}$,
every $\varepsilon>0$, and
every compact subset $K$ of $\R^n$, there exists a continuous
function $\rho_{\varepsilon,K}:\R^n\rightarrow {\R}$ satisfying
(\ref{eq12}) such that, for all
$e,d: \R^n \times [0,+\infty) \rightarrow \R^n$
satisfying the regularity assumptions (\ref{standing}) and
\begin{equation}
\begin{split}
& \sup_{[0,+\infty)}|e(x,\cdot)|\leq \min(\rho_{opt}(d(x,\mathcal{S})),
\rho_{\varepsilon,K}(x))
,\\
& {\mathrm{esssup}}_{[0,+\infty)}|d(x,\cdot)|\leq
\min(\rho_{opt}(d(x,\mathcal{S}),
\rho_{\varepsilon,K}(x))   ,\ \forall
x\in\R^n,
\end{split}
\end{equation}
and for every $x_0\in K \cap (\R^n\setminus \Om)$, the solution of
(\ref{ch:eq:5}), starting from
$x_0$, reaches $\bar{x}$ within time $T_{\bar{x}}(x_0)+\varepsilon$.
\end{itemize}
\end{lemma}

\begin{proof}
Since Carath\'eodory conditions hold for the system
(\ref{ch:eq:5}), the existence of a unique
Carath\'eodory solution of (\ref{ch:eq:5}), for every initial condition,
is ensured.
The inequality (\ref{eq6}) follows from the constraint (\ref{contrainte}).
Since the optimal controller $u_{opt}$ defined by
(\ref{eq20}) is continuous on $\R^n\setminus\mathcal{S}$,
Lemma \ref{leminvariance2} implies the existence of $\rho_{opt}:
\R^n\rightarrow [0,+\infty)$ so that the {\em robust stability} and the
{\em finite time convergence} properties hold.

The so-called {\em optimality} property follows from the definition of
$u_{opt}$, from the
continuity of solutions with respect to disturbances, and from the
compactness of the set of all solutions starting from $K\cup (\R^n\setminus
\Om)$.
\end{proof}


\subsection{Definition of the hybrid feedback law}\label{sechybrid}
A switching strategy must be defined in order to connect the first
component (optimal controller), and the second component (consisting a
a set of controllers, stated in Lemma \ref{e160}). The switching
strategy is achieved by adding a dynamical discrete variable $s_d$ and
using a hybrid feedback law, described next.

\subsubsection{Definitions}\label{secdefgeneralhybrid}
Let $\FF=\{1,\ldots,\5\}$, and $\NNN$ be a countable set. In the
sequel, greek letters refer to elements of $\NNN$. Fix
$\infini$ an element of $\NNN$. We emphasize that we do not introduce any
order in $\NNN$. However, intuitively, we consider that $\infini$ is
the {\em largest} element of $\NNN$, \textit{i.e.}, $\infini$ is 
{\em greater} than any other element of $\NNN$ (see in particular Remark
\ref{e164bis} below).

Given a set-valued map $F:\R^n\rightrightarrows \R^n$, we
define the solutions $x(\cdot)$ of the differential inclusion
$\dot x \in F(x) $
as all absolutely continuous functions satisfying $\dot x(t)\in
F(x(t))$ almost everywhere.
\begin{definition}
\label{def3} The family $(\R^n\setminus\{\bar x \},((\Omega_{\alpha,l})_{l \in
\FF},g_{\alpha})_{\alpha\in \NNN})$ is said to satisfy the
property $(\mathcal{P})$ if:
\begin{enumerate}
\item for every $(\al,l)\in \NNN\times \FF$, the set
  $\Omega_{\alpha,l}$ is an open subset of $\R^n$;
\item for every $\al \in\NNN$, and every $m>l\in \FF$,
  \begin{equation}
    \label{eq:16}
    \Om_{\al,l}\subsetneq\clo (\Om_{\al,l})\subsetneq \Om_{\al,m} ;
  \end{equation}
\item for every $\al$ in $\NNN$, $g_\al$ is a smooth vector field,
  defined in a neighborhood of $\clo (\Om_{\al,\5})$, taking values in
  $\R^n$;
\item
for every $(\al,l)\in \NNN\times\FF$, $l\leq \4$, there exists
 a continuous function $\rho_{\al,l}:\R^n \ra [0,+\infty)$ satisfying
 $\rho_{\al,l}(x)\neq 0 $ whenever $x\neq \bar x$ such that
every maximal solution $x(\cdot )$
 of
\begin{equation}
\label{e113}
\dot x \in  g_\al(x)+B(0,\rho_{\al,l}(x))  ;
\end{equation}
defined on $[0,T)$ and starting from $\partial \Om_{\al,l}$, is
such that
$$
x(t)\in\clo (\Om_{\al,l+1})  ,\  \fa t \in [0,T) ;
$$
\item for every $l\in \FF$, the sets $(\Om_{\al,l})_{\al \in \NNN}$ form a
locally finite covering of $\R^n\setminus\{\bar x \}$. 
\end{enumerate}
\end{definition}
\begin{remark}\label{e164}
Some observations are in order.
\begin{itemize}
\item First note that this notion is close to the notion of a family
  of nested patchy vector fields defined in
\cite{pat}. However note that, in general, the sets
$(\Om_{\alpha,l}, g_\alpha)$ may not be a patch as defined in
\cite{anbr99,pat}. Indeed, due to the property 4, the set
$\Omega_{\alpha,l}$ may not be invariant for the system
(\ref{e113}). Since the notion of a patch is one of
the main ingredients of the proofs of \cite{cdc05}, we cannot apply
\cite{cdc05} directly, even though some notions are however in
common (see in particular Definition \ref{def4b} below).
\item On the one hand, the function $\rho_{\al,l}$
allows to get robustness with respect to external disturbances. On
the other hand, the gap between the different patches given by
(\ref{eq:16}) allow to get robustness with respect to measurement
noise (see Definition \ref{e126} below for a precise statement of
an admissible radius of measurement noise and external
disturbances).
\item To state our main result, we need
consider a family of three nested patchy vector fields. The
patches 1, 2, 3, 4 and 6 define the dynamics of the discrete
component of our hybrid controller (see Definition \ref{def4b}
below). The patch 5 is used for technical reasons to handle the
measurement noise.
\end{itemize}
\end{remark}

We next define a class of
hybrid controllers as those considered in Section \ref{sec2} (see also
\cite{cdc05}).

\begin{definition}
\label{def4b} Let $(\R^n\setminus\{\bar x\},((\Omega_{\alpha,l})_{l \in
\FF},g_{\alpha})_{\alpha   \in \NNN})$ satisfy the property
$(\mathcal{P})$ as in Definition \ref{def3}. Assume that, for every
$\alpha$ in $\NNN$, there exists a smooth function $k_\alpha$ defined
in a neighborhood of $\Om_{\alpha,\5}$ and taking values in $\R^m$, such that,
for every $x$ in a neighborhood of $\Om_{\alpha,\5}$,
\begin{equation}
\label{e73b}
    g_\alpha(x)=f(x,k_\alpha(x))  .
\end{equation}
Set
\begin{eqnarray}
&D_{1}= \Om_{\infini,2}
  ,
\\
&D_{\alpha,2}=\R^n\setminus \Om_{\alpha,\4}  . \label{e185}
\end{eqnarray}
Let $(C,D,k,k_d)$ be the hybrid feedback defined by
\begin{eqnarray}
& \label{e140b} C= \Big\{(x,\alpha)\ \vert\ x\in
\Big(\clo(\Om_{\alpha,4})\setminus
\Omega_{\infini,1}\Big)  \Big\} ,
\\
&\label{e144b} D=\{ (x,\alpha)\ \vert\ x\in D_{1} \cup
D_{\alpha,2} \}  ,
\\
& \label{eq::4b}
\begin{array}{rclcl}
k:\R^n\times  \NNN &\ra& \R^m& \\
(x,\alpha)&\mapsto &k_\alpha (x)& \textrm{if} & x \in \Om_{\al,7},\\
&& 0 & \textrm{else},
\end{array}
\end{eqnarray}
and
\begin{equation}
\label{e10b}
\begin{array}{rcll}
k_d:\R^n\times \NNN&
\rightrightarrows& \NNN
\\
(x,\alpha)&\ms& \infini , & \mbox{if } x \in
\clo(\Om_{\infini,1} \cap D_{1})\mbox{ and if } x \not\in
D_{\alpha,2},
 \\
&&\alpha' , & \mbox{if } x \in \clo(\Om_{\alpha',1} \cap
D_{\alpha,2}).
\end{array}
\end{equation}

The 4-tuple $(C,D,k,k_d)$ is a hybrid feedback law on $\R^n$ as
considered in Section \ref{contr_sect}. We denote by
$\mathcal{H}_{(e,d)}$ the system (\ref{ch:e1}) in closed-loop with
such a feedback with the perturbations $e$ and $d$ as measurement
noise and external disturbance respectively.
\end{definition}

\begin{remark}\label{e164bis}
In this definition, we do not use any order in $\NNN$. However, in light of
\cite{cdc05}, we consider that $\infini$ is greater than any other element
of $\NNN$. This element $\infini$ has a particular role in the sequel,
since it will refer to the optimal controller in the hybrid feedback law.

This hybrid controller takes advantage of the existence of regions
where different controllers $k_\al$ exist and, roughly speaking,
allows the hybrid variable to choose between the different
controllers. This is the main idea of the hysteresis as done in
\cite{MCSS} to unit two controllers.
\end{remark}

Note that the concept of a hybrid feedback law of Definition
\ref{def4b} is similar to the one
of \cite{cdc05}. However, in \cite{cdc05}, the hybrid feedback
laws are derived from a family of patchy vector fields, whereas they
are here derived from a family
satisfying the property $(\mathcal{P})$ as considered in
Definition \ref{def3}.


\subsubsection{Properties of solutions}\label{secpropertiessolutionshybrid}
In this section, we investigate some properties of the solutions of the
system in closed-loop with the hybrid feedback law defined above.

\begin{definition}
\label{e126}
Let $\chi$ $:\R^n\ra\R$ be a continuous map such that $\chi(x)>0$, for every
$x\neq \bar x$.
\begin{itemize}
\item We say that $\chi$ is
an admissible radius for the measurement noise, if, for every
$x\in \R^n$ and every $\al\in \NNN$, such that $x\in
\Om_{\al,\5}$,
\begin{equation}
\label{eq::6} \chi(x)<\cfrac{1}{2} \min_{l \in
\{1,\ldots,\4\}}d(\R^n\setminus \Om_{\al,l+1},\Om_{\al,l})  .
\end{equation}
\item  We say that $\chi$ is an admissible radius for the external
  disturbances if, for every $x\in \R^n$, we have
$$
\chi(x)\leq \max_{(\al,l), \; x \in \Om_{\al,l}} \rho_{\al,l} (x) .
$$
\end{itemize}
\end{definition}

There exists an admissible radius for the measurement noise and
for the external disturbances (note that, from (\ref{eq:16}), the
right-hand side of the inequality (\ref{eq::6}) is positive).

Consider an admissible radius $\chi$ for the measurement noise and
the external disturbances. Let $e$ and $d$ be a measurement noise
and an external disturbance respectively, such that, for all
$(x,t)\in \R^n\times [0,+\infty)$,
\begin{equation}
\label{e138} e(x,t)\leq \chi(x) ,\ \ d(x,t)\leq \chi(x).
\end{equation}

The properties of the solutions of the system in
closed-loop with the hybrid feedback law defined in Definition
\ref{def4b} are similar to the ones of \cite{cdc05}. Hence, we skip
the proof of the following
three lemmas which do not use Statement 4 of Definition
\ref{def3}, but only the definition of the hybrid feedback law.
\begin{lemma}
\label{lem1}
For all $(x_0,s_0)\in\R^n\times\NNN$, there exists a
solution of $\mathcal{H}_{(e,d)}$ starting from $(x_0,s_0)$.
\end{lemma}

Recall that a Zeno solution is a complete solution
whose domain of
definition is bounded in the $t$-direction. A solution $(x,s_d)$, defined on
a hybrid domain $S$, is an instantaneous Zeno solution, if there exist
$t\geq 0$ and an infinite number of $j\in \N$ such that $(t,j)\in S$.

The Zeno solutions do not require a special treatment.
\begin{lemma}\label{e167b}
There do not exist
instantaneous Zeno solutions, although a finite number of switches may occur
at the same time.
\end{lemma}
%
%
%
%
%
%
%

We note, as usual, that maximal solutions of $\mathcal{H}_{(e,d)}$
blow up if their domain of definition is bounded.
Since Zeno solutions are avoided, the blow-up phenomenon
only concerns the $t$-direction of the domain of definition, and
we get the following result (see also \cite[Prop. 2.1]{go-te}).

\begin{lemma}
\label{th10}
Let $(x,s_d)$ be a maximal solution of $\mathcal{H}_{(e,d)}$ defined on
a hybrid time $S$. Suppose that the supremum $T$ of $S$ in the $t$-direction
is finite. Then,
$$
\limsup_{t\ra T, (t,l)\in S} |x(t,l)| =+\infty .
$$
\end{lemma}

We conclude this series of technical lemmas by studying the
behavior of solutions between two jumps. For every $\al\in \NNN$, set
\begin{equation}
\label{e82}
\begin{array}{rcl}
\tau_{\al}&=&\sup \big\{T\ \vert\ x \mbox{ is a Carath\'eodory solution of }
\dot x = f(x,k_{\al})+B(0,\chi(x))
\\
&&\qquad\qquad  \mbox{ with } x(t) \in \Om_{\al,\5}, \; \fa t\in
[0,T)\big\} \ .
\end{array}
\end{equation}
Note that there may exist $\al\in\NNN$ such that $\tau_\al=+\infty$.

\begin{lemma}
\label{e171} Let $(x,s_d)$ be a solution of $\mathcal{H}_{(e,d)}$
defined on a hybrid time domain $S$ and starting in $(\R^n\setminus \{\bar
x\})\times \NNN$. Let $T$ be the supremum in the $t$-direction of $S$. Then,
one of the two following cases may occur:
\begin{itemize}
\item either there exists no positive jump time, more precisely there exists $\alpha \in \NNN$
 such that,
\begin{enumerate}
\item\label{e172} for almost every
$t\in(0,T)$ and for every $l$ such that $(t,l)\in S$,
one has $k(s_{d}(t,l))=k_{\alpha}$;
\item\label{e173} the map $x$ is a Carath\'eodory
solution of $\dot x=f(x,k_{\alpha})+d$ on $(0,T,)$;
\item\label{e201} for every $t\in(0,T)$, and every $l$ such that $(t,l)\in
S$, one has $x(t,l)+e(x(t,l),t)\in
\clo(\Om_{\alpha,4})\setminus\Om_{\infini,1}$;
\item\label{e176ter} for all $(t,l)\in S$, $t>0$, one has
$x(t,l)+e(x(t,l),t)\not\in D$, where $D$ is defined by (\ref{e144b});
\item\label{e181} the inequality $T<\tau_{\alpha}$ holds.
\end{enumerate}
\item or there exists a unique positive jump time, more precisely there
exist $\alpha \in \NNN\setminus\{\infini\}$ 
and $t_1\in (0,T)$ such that, letting $t_0=0$, 
$t_2=T$, $\alpha_0=\alpha$, $\alpha_1=\infini$, for every $j=0 ,1$,
the following properties hold:
\begin{enumerate}
\item[6.]\label{e166} for almost every
$t\in(t_j,t_{j+1})$ and for every $l$ such that $(t,l)\in S$,
one has $k(s_{d}(t,l))=k_{\alpha_j}$;
\item[7.]\label{e167} the map $x$ is a Carath\'eodory
solution of $\dot x=f(x,k_{\alpha_j})+d$ on $(t_j,t_{j+1})$;
\item[8.]\label{e168} for every $t\in(t_0,t_{1})$,
and every $l$ such that $(t,l)\in
S$, one has $x(t,l)+e(x(t,l),t)\in \clo(\Om_{\alpha,4})\setminus
\Om_{\infini,1}$;
\item[9.]\label{e169} for every $t$ in $(t_j,t_{j+1})$, 
and every $l$ such that $(t,l)\in S$,
one has $x(t,l)+e(x(t,l),t)\not\in D_{\alpha_j,2}$, where
$D_{\alpha_j,2}$ is defined by (\ref{e185});
\item[10.]\label{e170} the inequality $t_1<\tau_{\alpha_j}$ holds.
\end{enumerate}
\end{itemize}
\end{lemma}

\begin{proof}
Consider the sequence $(t_j)_{j \in\overm}$ of jump times,
\textit{i.e.}, the times such that $t_0=0$ and, for every $j \in
\overm$, $j\leq m-1$,
\begin{eqnarray} &t_{j}\leq t_{j+1}\ ,
\\
\label{e184b}
&(x(t_{j+1},j)+e(x(t_{j+1},j),t_{j+1}),s_d(t_{j+1},j))\in D \ ,
\end{eqnarray}
and
\begin{equation}
\label{e149}
(x(t_{j+1},j+n_j)+e(x(t_{j+1},j),t_{j+1}),s_d(t_{j+1},j+n_j))\in C
\ ,
\end{equation}
where $n_j$ is the finite number of instantaneous switches (see
Lemma \ref{e167b}). Let $\sigma:\N \rightarrow \N$ be
an increasing function such that
$t_{\sigma(j)}<t_{\sigma(j+1)}$.

Between two jumps, $s_d(t)$ is
constant, and thus, there exists a sequence $(\alpha_j)$ in $\NNN$ such
that, for every $ t \in (t_{\sigma(j)},t_{\sigma(j+1)})$, except for a
finite number of $t$, we have
\begin{eqnarray}\label{e147}
&s_{d}(t,\sigma(j))=\alpha_j  ,
\\
\label{e142b} &\mbox{$x$ is a Carath\'eodory solution of $\dot
x=f(x,k_{\alpha_{j}})+d$ on $(t_{\sigma(j)},t_{\sigma(j+1)})$}  ,
\end{eqnarray}
and
\begin{equation}
\label{e141b} k(s_{d}(t,\sigma(j)))=k_{\alpha_{j}}  .
\end{equation}
 From (\ref{e140b}), (\ref{e149}) and (\ref{e147}), we have, for every
$ t\in [t_{\sigma(j)},t_{\sigma(j+1)}]$,
\begin{equation}
\label{e143b} x(t,j)+e(x(t,j),t)\in
\clo(\Om_{q_{j},4})\setminus\Om_{\infini,1}  .
\end{equation}
Note that, from (\ref{e141b}), (\ref{e143b}), and Statement 4
of Definition \ref{def3}, for every $t>0$ such that $t\in
[t_{\sigma(j)},t_{\sigma(j+1)}]$, one has 
\begin{equation}
\label{e184}
x(t,\sigma(j))\not\in
D_{\alpha_j,2}.
\end{equation}
Therefore, the positive jump time may occur only at time $t_j$ where
the point $x(t_j,l)+e(x(t_j,l),t_j)$ belongs to $D_1$. Thus, there
exists at most one positive
jump time.
 From (\ref{e142b}) and (\ref{e141b}), Statements
\ref{e172}, \ref{e173}, 6, 
7 
hold.
Statements \ref{e201} and 8 
are deduced from (\ref{e143b}).
Equation (\ref{e184}) implies Statements \ref{e176ter} and 9. 
Finally, Statements
\ref{e181} and 10 
are a consequence of (\ref{e82}) and (\ref{e142b}).
\end{proof}



\subsubsection{Definition of the hybrid feedback law, and switching
  strategy}\label{secswitching}
We next define our hybrid feedback law. Let $\eps>0$ and $K$
be a compact subset of $\R^n$. Let $\Omega$ be the neighborhood of
$\mathcal{S}$ given by Lemma \ref{e160}. For this neighborhood
$\Omega$, let $\Omega'\subset\Omega$ be the neighborhood of
$\mathcal{S}$ yielded by Lemma \ref{leminvariance2}.

Let $\mathcal{N}$  be the countable set defined by 
$$
\mathcal{N} = \{(i,p) , \ i \in \N , \ p \in \mathcal{N}_i\} \cup
\{\infini\}  ,
$$
where $\infini$ is an element of $\N\times \N$, distinct from every
$(i,p)$, $i\in\N$, $p\in\NNN_i$.

We proceed in two steps. We first define $k_\al$ and $\Om
_{\al,l}$, where $\al \in \NNN \setminus \{\infini\}$ and $l\in
\mathcal{F}$. Then, we define $k_\infini$ and $\Om _{\infini,l}$, where
$l\in \mathcal{F}$.

\begin{enumerate}
\item Let $i\in\N$.
Lemma \ref{e160}, applied with the stratum $M_i$,
implies the existence of a family of smooth
controllers $(k_{i,p})_{p\in\NNN_i}$
satisfying the constraint (\ref{contrainte}), and
of a family of neighborhoods $(\Om_{i,p,7})_{p\in\NNN_i}$.
The existence of the families $(\Om _{i,p,1})_{p\in\NNN_i},\ldots,(\Om
_{i,p,6})_{p\in\NNN_i}$, satisfying
\begin{equation*}
    \Om_{i,p,l}\subsetneq\clo (\Om_{i,p,l})\subsetneq \Om_{i,p,m} ,
  \end{equation*}
for every $m>l\in \FF$,
follows from a finite induction argument,
using Lemma \ref{e160}.


We have thus defined $k_{i,p}$ and $\Om _{i,p,l}$, where $(i,p) \in
\NNN \setminus \{\infini\}$ and $l\in \mathcal{F}$.

\begin{remark}
It follows from \cite{Ag} that, near the point $\bar{x}$, the cut
locus $\mathcal{S}$ is contained in a conic neighborhood
$\mathcal{C}$ centered at $\bar{x}$ (as shaped on Figure
\ref{fighybrid}), the axis of the cone being transversal to the
subspace $\mathrm{Span}\{f_1(\bar{x}),\ldots, f_m(\bar{x})\}$.
Hence, up to modifying slightly the previous construction, we
assume that, near $\bar{x}$, the set $\bigcup_{\al \in \NNN\setminus \{\infini\}, \; l \in \FF } \Om_{\al,l}$ is contained in this conic neighborhood.
\end{remark}

\medskip

\item
It remains to define the sets $\Om _{\infini,l}$, where $l\in
\mathcal{F}$, and the controller $k_\infini$. Let $\Om_{\infini,1}$ be an
open set of $\R^n$ containing $\R^n\setminus \bigcup _{\al\in\NNN\setminus
\{\infini\}}\Om_{\al,1}$ and contained in
$\R^n\setminus\mathcal{S}$. From the previous remark, the point
$\bar{x}$ belongs to $\clo(\Om_{\infini,1})$. Lemma \ref{leminvariance2},
applied with $\Om=\R^n\setminus \clo(\Om_{\infini,1})$, allows to
define $k_\infini$ as $k_{opt}$, and $\Om'$ a closed subset of $\R^n$
such that
\begin{equation}
\label{17juinbis}
\Omega'\subsetneq \Omega,
\end{equation}
and such that $\Om'$ is a neighborhood of $\mathcal{S}$. Set
$\Om_{\infini,2}=\R^n\setminus \Om'$; it is an open subset of
$\R^n$, contained in $\R^n\setminus \mathcal{S}$.
Moreover, from (\ref{17juinbis}),
$$
\Omega_{\infini,1}\subsetneq\clo (\Omega_{\infini,1})\subsetneq
\Omega_{\infini ,2}.
$$
The existence of the sets $\Om _{\infini,3},\ldots,\Om _{\infini,\5}$
follows from a finite induction argument, using Lemma
\ref{leminvariance2}. Moreover, from Lemma \ref{ch:le1}, we have
the following property: for every $l\in \{1, \ldots, 6\}$, for
every $x_0\in \Omega_{\infini,l}$, the unique Carath\'eodory solution
$x(\cdot )$ of (\ref{ch:eq:5}), with $x(0)=x_0$, satisfies
$x(t)\in \Om_{\infini,l+1}$, for every $t\geq 0$.
\end{enumerate}
\medskip
Therefore, $(\R^n\setminus\{\bar x\},((\Omega_{\alpha,l})_{l \in
\FF},g_{\alpha})_{\alpha   \in \NNN})$ is a family satisfying the
property $(\mathcal{P})$ as in Definition \ref{def3},
where $g_\alpha$ is a function defined in a neighborhood of
$\Om_{\alpha,\5}$ by
$$
    g_\alpha(x)=f(x,k_\alpha) .
$$
The hybrid feedback law $(C,D,k,k_d)$ is then defined according to Definition
\ref{def4b}.


\subsection{Proof of Theorem \ref{CH:TH1}}\label{e157}
Let $\eps>0$, and $K$ be a compact subset of $\R^n$.
Consider the hybrid feedback law $(C,D, k,k_d)$ defined
previously. Let $\chi$ be an admissible radius for the external
disturbances and the measurement noise (see Definition
\ref{e126}). Up to reduce this function, we assume that, for
every $\al\in\NNN\setminus \{\infini\}$,
\begin{eqnarray}
&\chi(x)\leq \rho_{opt}(d(x,\mathcal{S})), \quad \fa x \in
\Om_{\infini,\5},
\\
&\chi(x)\leq \rho _{\al}(x)  , \quad \fa x \in \Om_{\al,\5}.
\end{eqnarray}
Note that, from the choice of the components of the
hybrid feedback law, and from Lemmas \ref{e160} and \ref{e171},
for every
 $\al\in \NNN\setminus \{\infini\}$, the constant $\tau _{\al}$
 defined by (\ref{e82}) is such that $\tau _{\al} <\eps$.

Let us prove that the point $\bar{x}$ is a semi-global quasi-minimal time
robust
stable equilibrium for the system $\mathcal{H}_{(e,d)}$ in closed-loop with
the hybrid feedback law $(C,D,k,k_d)$
as stated in Theorem \ref{CH:TH1}.

\par\vspace{1em}\noindent{\em Step~1: Completeness and global stability}\\
Let $R>0$ and $\delta: [0,+\infty)\rightarrow [0,+\infty)$ of class
${\mathcal{K}}_\infty$ be such that, for every $\al\in\NNN\setminus
\{\infini\}$,
\begin{eqnarray}
&\delta(x)\leq \delta_{opt}(R) , \quad \fa x \in \Om_{\infini,\5},
\\
&\delta(x)\leq \delta _{\al}(R) , \quad \fa x \in \Om_{\al,\5} ,
\end{eqnarray}
where the functions $\delta _{\al}$ are defined in Lemma \ref{e160}.
Let $e$, $d$ be two functions satisfying the regularity assumptions and
(\ref{e138}). Let $(x,s_d)$ be a maximal solution of
$\mathcal{H}_{(e,d)}$ on a hybrid domain $S$
starting from $(x_0,s_0)$, with $|x_0|<R$.
\relax{}From Lemmas \ref{ch:le1}, \ref{e160} and \ref{e171}, we have,
for every $(t,l)\in S$,
\begin{equation}\label{bidule}
|x(t,l)-\bar{x}|\leq \delta(R) .
\end{equation}
Therefore, the conclusion of Lemma \ref{th10} cannot hold (since
$\limsup_{t\ra
T, (t,l)\in S}|x(t,l)|\neq +\infty$), and thus,
the supremum $T$ of $S$ in the $t$-direction
is infinite, and the maximality property follows.
The stability property follows from (\ref{bidule}).

\par\vspace{1em}\noindent{\em Step 2: Uniform finite time convergence
property}\\
Let $x_0\in B(\bar{x},R)$, and $s_0\in \NNN$. Let
$(x,s_d)$ denote the solution of
$\mathcal{H}_{(e,d)}$ starting from $(x_0,s_0)$.
If $x_0=\bar{x}$, then, using (\ref{eq::4b}) and the fact that
$\chi(\bar{x})=0$, the solution remains at the point $\bar{x}$.
We next assume that $x_0\neq \bar{x}$.
Let $\alpha_0\in \NNN$ such that
$x(\cdot )$ is a solution of $\dot x=f(x,k_{\al_0}(x))+d$ on
$(0, t_1)$ for a $t_1>0$ given by Lemma \ref{e171}.

If $\alpha_0=\infini$, then the feedback law under consideration
coincides with the optimal controller and, from Statement
\ref{e176ter} of Lemma \ref{e171}, there does not exist any switching
time $t>0$. Then, from Lemma \ref{ch:le1}, $x(\cdot )$ reaches
$\bar{x}$ within time $T_{\bar{x} }(x_0) +\eps$.

If $\alpha_0\neq \infini$, then, from Lemmas \ref{e160} and \ref{e171}, the
solution $x(\cdot )$ leaves $\Om_{\al_0,\5}$ within time $\eps$ and then
enters the set $\Om_{\infini,\5}$. Therefore, since $\tau _{\al} <\eps$,
$x(\cdot)$ reaches $\bar{x}$ within time $T_{\bar{x} }(x_1) +\eps$,
where $x_1$ denotes the point of $x(\cdot )$ when
entering $\Om_{\infini,\5}$.

Let $\tau(R)= \max _{x \in\delta(R)} T(x) + \eps$. With (\ref{bidule}), we
get (\ref{eq11}) and the uniform finite time property.
Note that, from Lemma \ref{e160}, the constraint (\ref{eq3}) is satisfied.

\par\vspace{1em}\noindent{\em Step 3: Quasi-optimality}\\
Let $K$ be a compact subset of $\R^n$, and $(x_0,s_0) \in K\times
\NNN$. Let $R>0$ such that $K\subset B(0,R)$. \relax{}From the
previous arguments, two cases occur:
\begin{itemize}
\item the solution starting from $(x_0,s_0)$
reaches
$\bar{x}$ within time $T_{\bar{x} }(x_0)+\eps$ whenever $\alpha_0=\infini$;
\item the solution starting from $(x_0,s_0)$
reaches $\bar{x}$ within time $T_{\bar{x} }(x_1) +\eps$, whenever
$\alpha_0\neq \infini$, where $x_1$ denotes the point of
$x(\cdot )$ when entering $\Om_{\infini,\5}$. Up to reducing the
neighborhoods $\Om_{\al,l}$, one has $|T_{\bar{x}
}(x_0)-T_{\bar{x} } (x_1) |\leq \eps$. Indeed, from Remark
\ref{remequivSR}, the function $T_{\bar{x}}(\cdot)$ is uniformly
continuous on the compact $K$.
\end{itemize}
Hence, the maximal solution starting from $(x_0,s_0)$
reaches $\bar{x}$ within time
$T_{\bar{x}}(x_0)+2\eps$. This is the quasi-optimality property.

\par\vspace{1em}
Theorem \ref{CH:TH1} is proved.

\section*{Acknowledgments} The authors are grateful to Ludovic
Rifford for constructive comments and suggestions.


\end{document}